\documentclass[12pt,reqno]{article}

\usepackage[usenames]{color}
\usepackage{amssymb}
\usepackage{amsmath}
\usepackage{amsthm}
\usepackage{amsfonts}
\usepackage{amscd}
\usepackage{graphicx}

\usepackage[colorlinks=true,
linkcolor=webgreen,
filecolor=webbrown,
citecolor=webgreen]{hyperref}

\definecolor{webgreen}{rgb}{0,.5,0}
\definecolor{webbrown}{rgb}{.6,0,0}

\usepackage{color}

\usepackage{graphics}
\usepackage{latexsym}

\begin{document}

\theoremstyle{plain}
\newtheorem{theorem}{Theorem}
\newtheorem{corollary}[theorem]{Corollary}
\newtheorem{lemma}[theorem]{Lemma}
\newtheorem{proposition}[theorem]{Proposition}

\theoremstyle{definition}
\newtheorem{definition}[theorem]{Definition}
\newtheorem{example}[theorem]{Example}
\newtheorem{conjecture}[theorem]{Conjecture}

\theoremstyle{remark}
\newtheorem{remark}[theorem]{Remark}

\begin{center}
\vskip 1cm{\LARGE\bf On new divisibility properties of generalized central trinomial coefficients and Legendre polynomials.
\vskip 1cm}
\large
Jovan Miki\'{c}\\
University of Banja Luka\\
Faculty of Technology\\
Bosnia and Herzegovina\\
\href{mailto:jovan.mikic@tf.unibl.org}{\tt jovan.mikic@tf.unibl.org} \\
\end{center}

\vskip .2in

\begin{abstract}
We present a new formula for the highest power of $a+b$ that divides the sum $B(n,m,a,b)=\sum_{k=0}^{n}\binom{n}{k}^m a^{n-k}b^k$ for the case $m=2$. By using this formula, we give complete 3-adic valuation for central Dellanoy numbers.
Also, we find the highest power of an odd integer $x$ that divides Legendre's polynomial $P_{n}(x)$. By using the same idea, generalized trinomial coefficients and generalized Motzkin numbers are treated. As a result,  we give complete 3-adic valuation for little  Schr\"{o}der numbers and restricted hexagonal numbers.
By using new class of binomial sums, we examine divisibility of $B(n,m,a,b)$  by powers of $a+b$ for $m >2$.
\end{abstract}


\noindent \textbf{2020} {\it \textbf{Mathematics Subject Classification}}:
05A10, 11B65.

\section{Introduction}\label{s:1}

Let us consider the following sum:

\begin{equation}\label{eq:1}
B(n,m,a,b)=\sum_{k=0}^{n}\binom{n}{k}^m a^{n-k}b^k\text{;}
\end{equation}
where $n$ and $m$ are natural numbers such that $m\geq 2$, and $a$, $b$ are integers.

Recently \cite[Section 2, Proposition 1, p.\ 4]{TerAlin}, it is shown that for $m=2$, the sum $B(n,2,a,b)$ represents the number of weighted free Dyck paths of semilength $n$ from $(0,0)$ to $(2n,0)$ with up steps $(1,1)$  and down steps $(1,-1)$; where the weight $a$ is assigned to up steps in the peaks, and the weight $b$ to other up steps. Also, the sum $B(n,2,a,b)$ represents the number of the same weighted free Dyck paths of semilength $n$, where the weight $a$ is assigned to even ascents and the weight $b$ to even descents.

The sum $B(n,2,a,b)$ is closely related to well-known Legendre polynomials $P_n(x)$. It is well-known \cite[Introduction]{VG2} that Legendre polynomials have various explicit representations. We refer readers to \cite[p.\ 1--3]{Koepf} for seven different definitions of these polynomials. For example,
\begin{equation}\label{eq:2}
P_n(x)=\frac{1}{2^n}\sum_{k=0}^{n}\binom{n}{k}^2(x+1)^k(x-1)^{n-k}\text{.}
\end{equation}

It is readily verified that the two following formulae are true:
\begin{align}
P_n(x)&=B(n,2,\frac{x-1}{2},\frac{x+1}{2})\text{,}\label{eq:3}\\
B(n,2,a,b)&=(a-b)^nP_n(\frac{a+b}{a-b})\text{, } a\neq b\label{eq:4}\text{.}
\end{align}

For $m=2$, $a=b=1$, the sum $B(n,m,a,b)$ is equal to the central binomial coefficient $\binom{2n}{n}$. It is well-known that the integer $\binom{2n}{n}$ counts all lattice paths from $(0,0)$ to $(n,n)$ by using steps $(1,0)$ and $(0,1)$. 
For $m=3$, $a=b=1$, the sum $B(n,m,a,b)$ is equal to the $n$-th Franel number $f_n=\sum_{k=0}^{n}\binom{n}{k}^3$.

For $m=2$, $a=1$, and $b=2$, the sum $B(n,m,a,b)$ is equal to the central Delannoy number $D_n$. It is well-known that  $D_n$ counts  all lattice paths from $(0,0)$ to $(n,n)$ by using steps $(1,0)$, $(0,1)$, and $(1,1)$. Such paths are called royal paths \cite[Introduction]{JSch}. Sulanke \cite{Sulanke} gave, in a recreational spirit, a collection of $29$ configurations counted by these numbers. Banderier and Schwer \cite{BS} gave additional information about origin and use of this number sequence. 
Central Delannoy numbers have a close relation to Legendre polynomials due to the fact that $D_n=P_n(3)$. See also \cite{GH, Law, LMoser}. 

Furthermore, central Delannoy numbers $D_n$ (see \cite[Section 6.4]{RS} and \cite{Sulanke 2})  satisfy  the second-order recurrence:
\begin{equation}\label{eq:5}
(n+2)D_{n+2}=3(2n+3)D_{n+1}-(n+1)D_n\text{;}
\end{equation}
where $n$ is a non-negative integer.

The large Schr\"{o}der numbers $S_n$ count all lattice paths from $(0,0)$ to $(n,n)$ by using steps $(1,0)$, $(0,1)$, and $(1,1)$ such that never go above the diagonal $y=x$. Large Schr\"{o}der numbers  bear the same relationship  to the central Delannoy numbers as the Catalan numbers do to the central binomial coefficients. Recall that the $n$-th Catalan number  $C_n=\frac{1}{n+1}\binom{2n}{n}$ counts all lattice paths from $(0,0)$ to $(n,n)$ by using steps $(1,0)$ and $(0,1)$ such that never go above the diagonal $y=x$.
Stanley \cite[Exercise 6.39]{RS} gave $11$ combinatorial objects counted by large Schr\"{o}der numbers $S_n$. See also \cite[Introduction]{JSch}.

It is known that:

\begin{equation}\label{eq:6}
S_n=\frac{1}{2}(-D_{n-1}+6D_n-D_{n+1})\text{;}
\end{equation}
where $n$ is a natural number.

The little Schr\"{o}der numbers $s_n$ count the number of plane trees with a given set of leaves, the number of ways of inserting parentheses into a sequence, and the number of ways of dissecting a convex polygon into smaller polygons by inserting diagonals. Let $n$ be a natural integer. Then
$s_n=\frac{1}{2}S_n$.

Furthermore, for $m=2$,  it is known that:
\begin{equation}\label{eq:7}
B(n,2,a,b)=\sum_{k=0}^{\lfloor \frac{n}{2}\rfloor}\binom{n}{k}\binom{n-k}{k}a^kb^k(a+b)^{n-2k}\text{.}
\end{equation}

For $m=3$, there is MacMahon's identity:
\begin{equation}
B(n,3,a,b)=\sum_{k=0}^{\lfloor \frac{n}{2}\rfloor}\binom{n}{2k}\binom{2k}{k}\binom{n+k}{k}a^kb^k(a+b)^{n-2k}\text{.}\label{eq:8}
\end{equation}
Recently,  new combinatorial proofs of Eqns.~(\ref{eq:7}) and (\ref{eq:8})  are presented, as well as combinatorial proofs of some known representations of Legendre polynomials based on weighted enumeration of free lattice paths \cite{TerAlin}.

Sun \cite{Sun2} introduced the following sum:
\begin{equation}\label{eq:8.1}
T_n(a,b)=\sum_{k=0}^{\lfloor \frac{n}{2}\rfloor}\binom{n}{k}\binom{n-k}{k}a^k\cdot b^{n-2k}\text{,}
\end{equation}
where $a$ and $b$ are integers.
The number $T_n(a,b)$ represents the number of weighted free Motzkin's paths of the length $n$ form $(0,0)$ to $(n,n)$ by using up steps $(1,1)$, down steps $(1,-1)$, and level steps $(1,0)$; where the weight $a$ is assigned to up steps and the weight $b$ to level steps.
The number $T_n(a,b)$ is also known \cite{Chen} as a generalized central trinomial coefficient. They can be defined as coefficients of $x^n$ in the expansion of $(x^2+bx+a)^n$.

By the Eq.~(\ref{eq:7}), it follows that 
\begin{equation}\label{eq:8.2}
B(n,2,a,b)=T_n(ab, a+b)\text{.}
\end{equation}
Also Sun\cite{Sun2} introduced the following sum:
\begin{equation}\label{eq:8.3}
M_n(a,b)=\sum_{k=0}^{\lfloor \frac{n}{2}\rfloor}\binom{n}{2k}C_k \cdot a^k\cdot b^{n-2k}\text{,}
\end{equation}
where $a$ and $b$ are integers. The number  $M_n(a,b)$ is known as a generalized Motzkin number. Particularly, it is known that:
\begin{equation}\label{eq:8.4}
s_{n+1}=M_n(2,3)\text{,}
\end{equation}
where $n$ is a non-negative integer.

Let us consider the following sum $B(n,m,1,-1)=\sum_{k=0}^{n}\binom{n}{k}^m(-1)^k$, where $n$ and $m$ are natural numbers such that $m\geq 2$. Obviously, 
\begin{equation}
B(2n-1,m,1,-1)=0\text{.}\label{eq:9}
\end{equation}

By setting $n:=2n$, $a:=1$, and $b:=-1$ in the Eq.~(\ref{eq:7}), it follows that 
\begin{equation}\label{eq:10}
B(2n,2,1,-1)=(-1)^n\binom{2n}{n}\text{.}\quad (\text{Kummer's formula})
\end{equation}

By setting $n:=2n$, $a:=1$, and $b:=-1$ in the Eq.~(\ref{eq:8}), it follows that 
\begin{equation}\label{eq:11}
B(2n,3,1,-1)=(-1)^n\binom{2n}{n}\binom{3n}{2n}\text{.}\quad (\text{Dixon's formula})
\end{equation}

In $ 1998$, Calkin  proved \cite[Thm.\ 1, p.\ 1]{NC} that the alternating binomial sum $B(2n, m,1,-1)$ is divisible by $\binom{2n}{n} $ for all non-negative integers $n$ and all positive integers $m$.   In 2007, Guo, Jouhet, and Zeng proved, among other things, two generalizations of Calkin's result \cite[Thm.\ 1.2, Thm.\ 1.3, p.\ 2]{VG}. In $2010$, Cao and Pan  \cite{HCHP} proved an extension of Calkin's result  \cite[Thm.\ 1]{NC}. In $2020$,  Calkin's result \cite[Thm.\ 1]{NC} was proved by using new class of binomial sums \cite[Section 5]{JM2}. See also \cite[Section 8, p.\ 14]{JM1}.

Let $a$ and $b$ be relatively prime integers such that $a+b\neq 0$. Note that:
\begin{equation}\label{eq:12}
B(n,m,a,b)\equiv a^n\cdot B(n,m,1,-1) \mod (a+b)\text{.}
\end{equation}

By Eqns.~(\ref{eq:9}) and (\ref{eq:12}), it follows that $B(2n-1,m,a,b)$ is always divisible by $a+b$. By Eqns.~(\ref{eq:10}) and (\ref{eq:12}), it follows that $B(2n,2,a,b)$ is divisible by $a+b$ if and only if $\binom{2n}{n}$ is divisible by $a+b$. Similarly, 
by Eqns.~(\ref{eq:11}) and (\ref{eq:12}), it follows that $B(2n,3,a,b)$ is divisible by $a+b$ if and only if  $\binom{2n}{n}\binom{3n}{n}$ is divisible by $a+b$.

Furthermore, by Calkin's result  \cite[Thm.\ 1, p.\ 1]{NC}, if $\binom{2n}{n}$ is divisible by $a+b$ then $B(2n,m,a,b)$ is divisible by $a+b$ for all natural numbers $n$ and $m$.

Our goal is to examine divisibility of $B(n,m,a,b)$ by  powers of $a+b$, for $m\geq 2$. For $m=2$, we shall find the highest power of $a+b$ such that divides the sum $B(n,2,a,b)$.

\section{Main Results}\label{s:2}

Let $\omega_{x}(y)$ denote the highest power of an integer $x$ that divides an integer $y$; where $x\neq 0$ and $x\neq \pm 1$. If an integer $x$ is a prime number, then we shall use a notation $v_p(y)$ instead of $\omega_{p}(y)$.
For calculating  $\omega_{x}(y)$, we use the following formula:
\begin{equation}
\omega_{x}(y)=\min (\lfloor \frac{v_p(y)}{v_p(x)} \rfloor: p \text{ is a prime such that divides } x ) \text{.}\label{eq:13}
\end{equation}

Our first main result is:

\begin{theorem}\label{t:1}
Let $n$  be a non-negative integer. Let $a$ and $b$ be integers such that $\gcd(a,b)=1$, $a+b\neq 0$, and $a+b\neq \pm 1$. For $m=2$, we have:

\begin{align}
\omega_{a+b}(B(2n,2,a,b))&=\omega_{a+b}(\binom{2n}{n})\text{,}\label{eq:14}\\
\omega_{a+b}(B(2n+1,2,a,b))&=1+\omega_{a+b}((2n+1)\binom{2n}{n})\text{.}\label{eq:15}
\end{align}
\end{theorem}

Our second main result is:

\begin{theorem}\label{t:2}
Let $n$ be a non-negative integer. Let $a$ and $b$ be integers such that $\gcd(a,b)=1$, $a+b\neq 0$, and $a+b\neq \pm 1$. Let $m$ be an arbitrary natural number greater than $2$. 
Then
\begin{align}
 \omega_{a+b}(B(2n,m,a,b))&\geq \omega_{a+b}(\binom{2n}{n})\text{,}\label{eq:16}\\
\omega_{a+b}(B(2n+1,m,a,b))&\geq 1+\omega_{a+b}((2n+1)\binom{2n}{n})\text{.}\label{eq:17}
\end{align}
\end{theorem}

We present several applications of Thms. ~\ref{t:1} and \ref{t:2}. 

The first application is for calculating the highest power of two that divides central binomial coefficient.
\begin{corollary}\label{cor:1}
Let $n$ be a non-negative integer.Then
\begin{align}
v_2(\binom{4n+2}{2n+1})&=1+v_2(\binom{2n}{n})\text{,}\label{eq:18}\\
v_2(\binom{4n}{2n})&=v_2(\binom{2n}{n})\text{.}\label{eq:19}
\end{align}
\end{corollary}

By using Eqns.~(\ref{eq:18}) and (\ref{eq:19}), it can be shown the well-known fact that $v_2(\binom{2n}{n})$ is equal to the number of ones in binary representation of the number $n$. 

The second application is for estimation for the highest power of two that divides Franel number $f_n$.
\begin{corollary}\label{cor:2}
Let $n$ be a non-negative integer.Then
\begin{align}
v_2(f_{2n+1})&\geq1+v_2(\binom{2n}{n})\text{,}\label{eq:20}\\
v_2(f_{2n})&\geq v_2(\binom{2n}{n})\text{.}\label{eq:21}
\end{align}
\end{corollary}

The third application is for calculating the highest power of three that divides central Delannoy number $D_n$.
\begin{theorem}\label{t:3}
Let $n$ be a non-negative integer.Then
\begin{align}
v_3(D_{2n+1})&=1+v_3 (2n+1)+v_3(\binom{2n}{n})\text{,}\label{eq:22}\\
v_3(D_{2n})&=v_3(\binom{2n}{n})\text{.}\label{eq:23}
\end{align}
\end{theorem}

The fourth application is for calculating the highest power of three that divides large Schr\"{o}der number $S_n$.
\begin{theorem}\label{t:4}
Let $n$ be a non-negative integer.Then
\begin{align}
v_3(S_{2n+1})&=v_3(C_n)\text{,}\label{eq:24}\\
v_3(S_{2n+2})&=1+v_3(2n+1)+v_3(C_n)\text{.}\label{eq:25}
\end{align}
\end{theorem}

Note that $v_3(s_n)=v_3(S_n)$.

Recently, Lengyel, among other, gave the formula \cite[Theorem 10, p.\ 8]{TLeng} for $3$-adic valuation of central Dellanoy numbers. Note that it is asserted that the formula \cite[Theorem 10, p.\ 8]{TLeng} is true for any $n$ sufficiently large.
Also, Lengyel gave the formula \cite[Theorem 17, p.\ 19]{TLeng} for $3$-adic valuation of large Schr\"{o}der numbers $S_n$ for the case $n\equiv 1(\mod 3)$.

The fifth application is for calculating the highest power of an odd integer $x$ that divides Legendre polynomial $P_n(x)$.
\begin{corollary}\label{cor:3}
Let $n$ be a non-negative integer, and let $x$ be an odd integer such that  $x \neq \pm 1$.
Then
\begin{align}
\omega_x(P_{2n+1}(x))&=1+\omega_x((2n+1)\binom{2n}{n})\text{,}\label{eq:27}\\
\omega_x(P_{2n}(x))&=\omega_x(\binom{2n}{n})\text{.}\label{eq:28}
\end{align}
\end{corollary}

Furthermore, from our proof of Theorem \ref{t:1}, we can calculate the highest power of an integer $b$ that divides generalized central trinomial coefficient $T_n(a,b)$.

\begin{theorem}\label{t:5}

Let $n$  be a non-negative integer. Let $a$ and $b$ be integers such that $\gcd(a,b)=1$, $b \neq \pm 1$, and $b \neq 0$. Then:

\begin{align}
\omega_{b}(T_{2n}(a,b))&=\omega_{b}(\binom{2n}{n})\text{,}\label{eq:28.1}\\
\omega_{b}(T_{2n+1}(a,b))&=1+\omega_{b}((2n+1)\binom{2n}{n})\text{.}\label{eq:28.2}
\end{align}
\end{theorem}

Also, by using the same idea from our proof of Theorem \ref{t:1}, we can calculate the highest power of an integer $b$ that divides generalized Motzkin number $M_n(a,b)$.

\begin{theorem}\label{t:6}
Let $n$  be a non-negative integer. Let $a$ and $b$ be integers such that $\gcd(a,b)=1$, $b\neq 0$, and $b\neq \pm 1$. Then:

\begin{align}
\omega_{b}(M_{2n}(a,b))&=\omega_{b}(C_n)\text{,}\label{eq:28.3}\\
\omega_{b}(M_{2n+1}(a,b))&=1+\omega_{b}((2n+1)C_n)\text{.}\label{eq:28.4}
\end{align}
\end{theorem}

For a proof of Theorem \ref{t:1}, we use the Eq.~(\ref{eq:7}) and the following lemma:

\begin{lemma}\label{l:1}
Let $n$ be a natural number, and let $p$ be a prime number. Then
\begin{equation}\notag
v_p((2n+1)\binom{2n}{n})\leq n\text{.}
\end{equation}
\end{lemma}

Furthermore, for a proof of Theorem \ref{t:2}, we use Lemma \ref{l:1} and a new class of binomial sums \cite[Eqns.~(27) and (28)]{JM2} that we call $M$ sums.

\section{New Class of Binomials Sums}\label{s:3}

\begin{definition}\label{def:1}
Let $n$ and $m$ be  positive integers, and let $a$ and $b$ be integers. Let $S(n,m,a,b)=\sum_{k=0}^{n}\binom{n}{k}^m F(n,k,a,b)$, where $F(n,k,a,b)$ is an integer-valued function. Then  $M$ sums for the sum $S(n,m,a,b)$ are, as follows:
\begin{equation}\label{eq:29}
M_S(n,j,t;a,b)=\binom{n-j}{j}\sum_{v=0}^{n-2j}\binom{n-2j}{v}\binom{n}{j+v}^t F(n, j+v, a,b)\text{;}
\end{equation}
where $j$ and $t$ are non-negative integers such that $j\leq \lfloor \frac{n}{2} \rfloor$.
\end{definition}

Obviously, the following equation \cite[Eq.~(29)]{JM2},
\begin{equation}\label{eq:30}
S(n,m,a,b)=M_S(n,0,m-1;a,b)
\end{equation}
holds.

By setting $t:=0$ in the Eq.~(\ref{eq:29}), it follows that
\begin{equation}\label{eq:31}
M_S(n,j,0;a,b)=\binom{n-j}{j}\sum_{v=0}^{n-2j}\binom{n-2j}{v}F(n, j+v, a,b)\text{.}
\end{equation}

Let  $n$, $j$, $t$, $a$, and $b$ be same as in Definition \ref{def:1}.

It is known \cite[Th.\ 8]{JM2} that  $M$ sums satisfy the following recurrence relation:
\begin{equation}\label{eq:32}
M_S(n,j,t+1;a,b)=\binom{n}{j}\sum_{u=0}^{\lfloor\frac{n-2j}{2}\rfloor}\binom{n-j}{u}M_S(n,j+u,t;a,b)\text{.}
\end{equation}

In one particular situation, Relation \ref{eq:32} has a simple consequence \cite[Section 4]{JM2} which is important for us.

\begin{corollary}\label{cor:6}
Let $n$ be a fixed natural number, and let $a$,$b$ be fixed integers.
Let us suppose that, for some non-negative integer $t'$, $M_S(n,j,t';a,b)$ is divisible by an integer $q(n,a,b)$ for all $0\leq j \leq \lfloor\frac{n}{2}\rfloor$.
Then $M_S(n,j,t;a,b)$ is divisible by $q(n,a,b)$ for all $t\geq t'$ and for all $0\leq j \leq \lfloor\frac{n}{2}\rfloor$. Furthermore, $S(n,m,a,b)$ is divisible by $q(n,a,b)$ for all $m\geq t'+1$.
\end{corollary}

The idea is to calculate  $M$ sums for the sum $B(n,m,a,b)$ by using Eq.~(\ref{eq:31}) and Relation (\ref{eq:32}). 
By using binomial theorem, we show that:
\begin{equation}\label{eq:33}
M_B(n,j,0;a,b)=\binom{n-j}{j}(ab)^j(a+b)^{n-2j}\text{.}
\end{equation}

Furthermore, we prove a slight generalization of the Eq.~(\ref{eq:7}):
\begin{equation}\label{eq:34}
M_B(n,j,1;a,b)=\binom{n}{j}\sum_{k=0}^{\lfloor \frac{n-2j}{2}\rfloor}\binom{n-j}{k}\binom{n-j-k}{j+k}(ab)^{j+k}(a+b)^{n-2j-2k}\text{.}
\end{equation}

By setting $j:=0$ in the Eq.~(\ref{eq:34}), we obtain the Eq.~(\ref{eq:7}).

$M$ sums give an elementary proof of Calkin's result  \cite[Thm.\ 1, p.\ 1]{NC}.  It is known \cite[Eqns.~(22) and (25)]{JM2} that: 
\begin{align*}
M_{B}(2n,j,0;1,-1)&=\begin{cases}0, & \text{if } 0 \leq j <n;\\
(-1)^n, & \text{if } j=n.
\end{cases}\\
M_{B}(2n,j,1;1,-1)&=(-1)^n\binom{2n}{n}\binom{n}{j}\text{;}
\end{align*}
where $0\leq j \leq n$. It follows that $q_1(2n,1,-1)=\binom{2n}{n}$ and $t'=1$.

By using $M$ sums, we can derive a slight generalization of, the MacMahon's identity, Eq.~(\ref{eq:8}). Also we can derive formulas for $B(n,m,a,b)$ for the cases $m=4$ and $m=5$.

Note that there are several other applications (\cite[Background]{JM5},\cite{JM2}, and \cite{JM6})  of $M$ sums. 

The rest of the paper is structured as follows. 
In Section \ref{s:4}, we give a proof of Lemma \ref{l:1} by using Legendre's theorem. In Section \ref{s:5}, we give a proof of Theorem \ref{t:1}. In Section \ref{s:6}, we prove Theorem \ref{t:2} by using new class of binomials sums.
In Section \ref{s:7}, we give a proof of Corollary \ref{cor:1}. In Section \ref{s:8}, we give a proof of Corollary \ref{cor:2}. In Section \ref{s:9}, we prove Theorem \ref{t:3}. In Section \ref{s:10}, we give a sketch of proofs
of Theorem \ref{t:5} and Corollary \ref{cor:3}. In Section \ref{s:11}, we give a proof of Theorem \ref{t:6}. In Section \ref{s:12}, we give a proof of Theorem \ref{t:4} by using Theorem \ref{t:6}. In Section \ref{s:13}, we give concluding remarks.

\section{A Proof of Lemma \ref{l:1}}\label{s:4}

\textbf{The first case}:

 Let $p=2$. Let $(2n)!!$ denote $\prod_{k=1}^{n}(2k)$, and let  $(2n-1)!!$ denote $\prod_{k=1}^{n}(2k-1)$.

We have gradually:

\begin{align*}
v_2((2n+1)\binom{2n}{n})&=v_2(2n+1)+v_2(\binom{2n}{n})\\
&=v_2(\binom{2n}{n})\\
&=v_2(\frac{(2n)!}{(n!)^2})\\
&=v_2(\frac{(2n)!!\cdot (2n-1)!!}{(n!)^2})\\
&=v_2(\frac{2^n\cdot n!\cdot (2n-1)!!}{(n!)^2})\\
&=v_2(\frac{2^n\cdot (2n-1)!!}{n!})\\
&=v_2(2^n)+v_2((2n-1)!!)-v_2(n!)\text{.}
\end{align*}

From the last equation above, it follows that
\begin{equation}\label{eq:39}
v_2((2n+1)\binom{2n}{n})=n-v_2(n!)\text{.}
\end{equation}

By the Eq.~(\ref{eq:39}), it follows that $v_2((2n+1)\binom{2n}{n})\leq n$.
This proves the first case.

\textbf{The second case}:

Let $p$ be an odd prime. By Legendre's formula, it follows that:

\begin{equation}\label{eq:40}
v_p(n!)=\sum_{k=1}^{\infty}\lfloor\frac{n}{p^k} \rfloor
\end{equation}

Let $n=\sum_{l=0}^{k}a_l\cdot p^l$ be the expansion of $n$ in the base $p$, where $a_l$ are non-negative integers such that $a_l<p$ for all $ l=1,k$.
Let $s_p(n)$ denote the sum $\sum_{l=0}^{k}a_l$\text{.}

By Legendre's formula, it follows Legendre's theorem \cite[Theorem 1.1]{Mihet}:

\begin{equation}\label{eq:41}
v_p(n!)=\frac{n-s_p(n)}{p-1}\text{.}
\end{equation}

Let us consider the number $c(n,p)=\frac{(pn)!}{(n!)^p}$. Obviously, $c(n,p)$ is always an integer due to the fact that $c(n,p)$ is a central multinomial coefficient.

It is known that:

\begin{equation}\label{eq:42}
c(n,p)=\prod_{k=2}^{p}\binom{kn}{n}\text{.}
\end{equation}

By Legendre's theorem \ref{eq:41}, it follows that:

\begin{align}
v_p(c(n,p))&=v_p(\frac{(pn)!}{(n!)^p})\notag\\
&=v_p((pn)!)-pv_p(n!)\notag\\
&=\frac{pn-s_p(pn)}{p-1}-p\cdot \frac{n-s_p(n)}{p-1}\notag\\
&=\frac{p\cdot s_p(n)-s_p(pn)}{p-1}\text{.}\label{eq:43}
\end{align}

By using the fact $s_p(pn)=s_p(n)$, from the Eq.~(\ref{eq:43}), it follows that:

\begin{equation}\label{eq:44}
v_p(c(n,p))=s_p(n)\text{.}
\end{equation}

From the expansion of $n$ in the base $p$, using the obvious ineqality $p^l \geq 1$, it follows that $s_p(n)\leq n$. By using the Eq.~(\ref{eq:44}), we obtain that:
\begin{equation}\label{eq:45}
v_p(c(n,p))\leq n\text{.}
\end{equation}

Let $F(n,k)=\frac{1}{(k-1)n+1}\binom{kn}{n}$ denote generalized Catalan numbers or Fuss-Catalan numbers, where $n$ and $k$ are natural numbers. It is well-known that Fuss-Catalan number is always an integer \cite[Eq.~(17.1), p.\ 375]{Koshy}.  
Obviously, $F(n,2)=C_n$ and $F(n,3)=\frac{1}{2n+1}\binom{3n}{n}$. See also \cite{JM6}.

Now we use assumption that $p$ is an odd prime. By the Eq.~(\ref{eq:42}), it follows that $c(n,p)$ is divisible by $\binom{3n}{n}\binom{2n}{n}$. Note that $\binom{3n}{n}\binom{2n}{n}=F(n,3)\cdot (2n+1)\binom{2n}{n}$.
It follows that $c(n,p)$ is divisible by $(2n+1)\binom{2n}{n}$. Therefore,

\begin{equation}\label{eq:46}
v_p((2n+1)\binom{2n}{n})\leq v_p(c(n,p))\text{.}
\end{equation}

Eqns.~(\ref{eq:45}) and (\ref{eq:46}) prove the second case of Lemma \ref{l:1}.
This completes the proof of Lemma \ref{l:1}.

\begin{remark}\label{r:1}
Let $p$ be an odd prime. Then, by using the same idea from the second case of the proof of Lemma \ref{l:1}, it can be shown that:
\begin{equation}\label{eq:47}
v_p((\prod_{k=2}^{p-1}(kn+1))\binom{2n}{n})\leq n\text{.}
\end{equation}
\end{remark}

\begin{remark}\label{r:2}
The Eq.~(\ref{eq:44}) is true for all prime numbers $p$. Particulary, for $p=2$, it follows that $v_2(\binom{2n}{n})$ is equal to number of ones in the binary representation of the number $n$.
\end{remark}

\section{A Proof of Theorem \ref{t:1}}\label{s:5}

Firstly, we give a proof of the Eq.~(\ref{eq:14}).

\subsection{A Proof of the Eq.~(\ref{eq:14})}\label{s:5.1}

By setting $n:=2n$ in the Eq.~(\ref{eq:7}), we obtain that:

\begin{equation}\label{eq:48}
B(2n,2,a,b)=\sum_{k=0}^{n}\binom{2n}{k}\binom{2n-k}{k}\cdot (ab)^k\cdot (a+b)^{2n-2k}\text{.}
\end{equation}

By using the substitution $k=n-s$, the Eq.~(\ref{eq:48}) becomes gradually:

\begin{align}
B(2n,2,a,b)&=\sum_{s=0}^{n}\binom{2n}{n-s}\binom{n+s}{n-s}(a+b)^{2s} (ab)^{n-s}\text{,}\notag\\
&=\sum_{s=0}^{n}\binom{2n}{n+s}\binom{n+s}{2s} (a+b)^{2s}\cdot (ab)^{n-s}\text{.}\label{eq:49}
\end{align}

It is readily verified that \cite[Eq.~(1.4), p.\ 5]{Koshy}:

\begin{equation}\label{eq:50}
\binom{2n}{n+s}\binom{n+s}{2s}=\frac{\binom{2n}{n}\binom{n}{s}^2}{\binom{2s}{s}}\text{.}
\end{equation}

By using the Eq.~(\ref{eq:50}), the Eq.~(\ref{eq:49}) becomes, as follows:

\begin{equation}\label{eq:51}
B(2n,2,a,b)=\sum_{s=0}^{n}\bigl{(} \frac{\binom{2n}{n}\binom{n}{s}^2}{\binom{2s}{s}} \bigr{)} (a+b)^{2s}\cdot (ab)^{n-s}\text{.}
\end{equation}

Let $p$ be a prime number such that divides $a+b$.

For $s=0$, the first summand in the Eq.~(\ref{eq:51}) is $\binom{2n}{n}\cdot (ab)^n$. Since $a+b$ and $ab$ are relatively prime, 
it follows that
\begin{equation}\label{eq:52}
v_p(\binom{2n}{n}\cdot (ab)^n)=v_p(\binom{2n}{n}).
\end{equation}

Let $1\leq s \leq n$. Note that by Lemma \ref{l:1}, it follows that
\begin{equation}\label{eq:53}
v_p(\binom{2s}{s})\leq s\text{.}
\end{equation}

Then, for $s+1$-th summand the Eq.~(\ref{eq:51}), we obtain that:

\begin{equation}\label{eq:54}
v_p(\bigl{(} \frac{\binom{2n}{n}\binom{n}{s}^2}{\binom{2s}{s}} \bigr{)} (a+b)^{2s}\cdot (ab)^{n-s} )=v_p(\binom{2n}{n})+2v_p(\binom{n}{s})-v_p(\binom{2s}{s})+(2s)\cdot v_p(a+b)\text{.}
\end{equation}

By using the Ineq.~(\ref{eq:53}), from the Eq.~(\ref{eq:54}), we obtain that:

\begin{equation}\label{eq:55}
v_p(\bigl{(} \frac{\binom{2n}{n}\binom{n}{s}^2}{\binom{2s}{s}} \bigr{)} (a+b)^{2s}\cdot (ab)^{n-s} )\geq v_p(\binom{2n}{n})+s\cdot(2v_p(a+b)-1)\text{.}
\end{equation}

Obviously, for $1\leq s \leq n$, 

\begin{equation}\label{eq:56}
v_p(\binom{2n}{n})+s\cdot(2v_p(a+b)-1)\geq v_p(\binom{2n}{n})+s \geq  v_p(\binom{2n}{n})+1 \text{.}
\end{equation}

By Ineq.~(\ref{eq:55}) and (\ref{eq:56}), for $1\leq s \leq n$, it follows that:
\begin{equation}\label{eq:57}
v_p(\bigl{(} \frac{\binom{2n}{n}\binom{n}{s}^2}{\binom{2s}{s}} \bigr{)} (a+b)^{2s}\cdot (ab)^{n-s} )\geq  v_p(\binom{2n}{n})+1\text{.}
\end{equation}

By using Eqns. (\ref{eq:52}) and (\ref{eq:57}), we obtain that:

\begin{equation}\label{eq:58}
v_p(B(2n,2,a,b))=v_p(\binom{2n}{n})\text{,}
\end{equation}
where $p$ is a prime number that divides $a+b$.

Now, by using Relation \ref{eq:13}, it follows that:

\begin{align}
\omega_{a+b}(B(2n,2,a,b))&=\min (\lfloor \frac{v_p((B(2n,2,a,b))}{v_p(a+b)} \rfloor: p \text{ is a prime such that divides } a+b )\label{eq:59}\\
&=\min (\lfloor \frac{v_p(\binom{2n}{n})}{v_p(a+b)} \rfloor: p \text{ is a prime such that divides } a+b )\label{eq:60}\\
&=\omega_{a+b}{(\binom{2n}{n})}\text{.}\label{eq:61}
\end{align}

The Eq.~(\ref{eq:61}) completes the proof of the Eq.~(\ref{eq:14}).

\subsection{A Proof of the Eq.~(\ref{eq:15})}\label{s:5.2}

By setting $n:=2n+1$ in the Eq.~(\ref{eq:7}), we obtain that:

\begin{align}\
B(2n+1,2,a,b)&=\sum_{k=0}^{n}\binom{2n+1}{k}\binom{2n+1-k}{k}\cdot (ab)^k\cdot (a+b)^{2n+1-2k}\text{,}\notag\\
&=(a+b)\sum_{k=0}^{n}\binom{2n+1}{k}\binom{2n+1-k}{k}\cdot (ab)^k\cdot (a+b)^{2n-2k}\text{.}\label{eq:63}
\end{align}

By using the substitution $k=n-s$, the Eq.~(\ref{eq:63}) becomes gradually:

\begin{align}
B(2n+1,2,a,b)&=(a+b)\sum_{s=0}^{n}\binom{2n+1}{n-s}\binom{n+s+1}{n-s}(a+b)^{2s} (ab)^{n-s}\text{,}\notag\\
&=(a+b)\sum_{s=0}^{n}\binom{2n+1}{n+s+1}\binom{n+s+1}{2s+1} (a+b)^{2s}\cdot (ab)^{n-s}\text{.}\label{eq:64}
\end{align}

It is readily verified that  \cite[Eq.~(1.4), p.\ 5]{Koshy}

\begin{equation}\label{eq:65}
\binom{2n+1}{n+s+1}\binom{n+s+1}{2s+1} =\frac{(2n+1)\binom{2n}{n}\binom{n}{s}^2}{(2s+1)\binom{2s}{s}}\text{.}
\end{equation}

By using the Eq.~(\ref{eq:65}), the Eq.~(\ref{eq:64}) becomes, as follows:

\begin{equation}\label{eq:66}
B(2n+1,2,a,b)=(a+b)\sum_{s=0}^{n}\bigl{[} \frac{(2n+1)\binom{2n}{n}\binom{n}{s}^2}{(2s+1)\binom{2s}{s}} \bigr{]} (a+b)^{2s}\cdot (ab)^{n-s}\text{.}
\end{equation}

By using the Eq.~(\ref{eq:66}), it follows that

\begin{equation}\label{eq:67}
\omega_{a+b}(B(2n+1,2,a,b))=1+\omega_{a+b}(\sum_{s=0}^{n}\bigl{[} \frac{(2n+1)\binom{2n}{n}\binom{n}{s}^2}{(2s+1)\binom{2s}{s}} \bigr{]} (a+b)^{2s}\cdot (ab)^{n-s})\text{.}
\end{equation}

Let $P(2n+1,2,a,b)$ denote $\sum_{s=0}^{n}\bigl{[} \frac{(2n+1)\binom{2n}{n}\binom{n}{s}^2}{(2s+1)\binom{2s}{s}} \bigr{]} (a+b)^{2s}\cdot (ab)^{n-s}$.

Then, the Eq.~(\ref{eq:67}) can be rewritten as:

\begin{equation}\label{eq:68}
\omega_{a+b}(B(2n+1,2,a,b))=1+\omega_{a+b}(P(2n+1,2,a,b))\text{.}
\end{equation}

Let $p$ be a prime number such that divides $a+b$.

For $s=0$, the first summand in the sum $P(2n+1,2,a,b)$ is $(2n+1)\binom{2n}{n}\cdot (ab)^n$. Since $a+b$ and $ab$ are relatively prime, 
it follows that
\begin{equation}\label{eq:69}
v_p((2n+1)\binom{2n}{n}\cdot (ab)^n)=v_p((2n+1)\binom{2n}{n}).
\end{equation}

Let $1\leq s \leq n$. 

Then, for $s+1$-th summand the sum $P(2n+1,2,a,b)$ , it follows that

$v_p(\bigl{(} \frac{(2n+1)\binom{2n}{n}\binom{n}{s}^2}{(2s+1)\binom{2s}{s}}\bigr{)} (a+b)^{2s}\cdot (ab)^{n-s})$ is equal to 

\begin{equation}\label{eq:70}
v_p((2n+1)\binom{2n}{n})+2v_p(\binom{n}{s})-v_p((2s+1)\binom{2s}{s})+(2s)\cdot v_p(a+b)\text{.}
\end{equation}

By using Lemma \ref{l:1}, from the Eq.~(\ref{eq:70}), it follows that:

\begin{equation}\label{eq:71}
v_p(\bigl{(} \frac{(2n+1)\binom{2n}{n}\binom{n}{s}^2}{(2s+1)\binom{2s}{s}} \bigr{)} (a+b)^{2s}\cdot (ab)^{n-s} )\geq v_p((2n+1)\binom{2n}{n})+s\cdot(2v_p(a+b)-1)\text{.}
\end{equation}

Obviously, for $1\leq s \leq n$, 

\begin{equation}\label{eq:72}
v_p((2n+1)\binom{2n}{n})+s\cdot(2v_p(a+b)-1)\geq v_p((2n+1)\binom{2n}{n})+s \geq  v_p((2n+1)\binom{2n}{n})+1 \text{.}
\end{equation}

By Ineq.~(\ref{eq:71}) and (\ref{eq:72}), for $1\leq s \leq n$, it follows that:
\begin{equation}\label{eq:73}
v_p(\bigl{(} \frac{(2n+1)\binom{2n}{n}\binom{n}{s}^2}{(2s+1)\binom{2s}{s}} \bigr{)} (a+b)^{2s}\cdot (ab)^{n-s} )\geq  v_p((2n+1)\binom{2n}{n})+1\text{.}
\end{equation}

By using Eqns. (\ref{eq:69}) and (\ref{eq:73}), we obtain that:

\begin{equation}\label{eq:74}
v_p(P(2n+1,2,a,b))=v_p((2n+1)\binom{2n}{n})\text{,}
\end{equation}
where $p$ is a prime number that divides $a+b$.

Now, by using Relation \ref{eq:13}, it follows that:

\begin{align}
\omega_{a+b}(P(2n+1,2,a,b))&=\min (\lfloor \frac{v_p((P(2n+1,2,a,b))}{v_p(a+b)} \rfloor: p \text{ is a prime  divides } a+b )\label{eq:75}\\
&=\min (\lfloor \frac{v_p((2n+1)\binom{2n}{n})}{v_p(a+b)} \rfloor: p \text{ is a prime  divides } a+b )\label{eq:76}\\
&=\omega_{a+b}{((2n+1)\binom{2n}{n}})\text{.}\label{eq:77}
\end{align}

By using Eqns.~(\ref{eq:68}) and (\ref{eq:77}),  the Eq.~(\ref{eq:15}) follows. This completes the proof of Theorem \ref{t:1}.

\section{A Proof of Theorem \ref{t:2}}\label{s:6}

Firstly, we give a proof of the Eq.~(\ref{eq:33}).

\subsection{A Proof of the Eq.~(\ref{eq:33})}

Obviously, the sum $B(n,m,a,b)$ is an instance of the sum $S(n,m,a,b)$, where $F(n,k,a,b)=a^{n-k}b^k$.

By the Eq.~(\ref{eq:31}), it follows that
\begin{align*}
M_B(n,j,0;a,b)&=\binom{n-j}{j}\sum_{v=0}^{n-2j}\binom{n-2j}{v}a^{n-j-v}b^{j+v}\text{,}\\
&=\binom{n-j}{j}a^j\cdot b^j\sum_{v=0}^{n-2j}\binom{n-2j}{v}a^{n-2j-v}b^{v}\text{,}\\
&=\binom{n-j}{j}(ab)^j(a+b)^{n-2j}\text{,}&&(\text{By the binomial theorem})
\end{align*}
where $0\leq j \leq \lfloor \frac{n}{2} \rfloor$.
This proves the  Eq.~(\ref{eq:33}).

Secondly, we give a proof of the Eq.~(\ref{eq:34}).

\subsection{A Proof of the Eq.~(\ref{eq:34})}

By setting $t:=0$ and $S:=B$ in the Eq.~(\ref{eq:32}), it follows that
\begin{equation}\label{eq:80}
M_B(n,j,1;a,b)=\binom{n}{j}\sum_{u=0}^{\lfloor\frac{n-2j}{2}\rfloor}\binom{n-j}{u}M_B(n,j+u,0;a,b)\text{.}
\end{equation}

By using the Eq.~(\ref{eq:33}), the Eq.~(\ref{eq:80}) becomes:
\begin{equation}\label{eq:81}
M_B(n,j,1;a,b)=\binom{n}{j}\sum_{u=0}^{\lfloor\frac{n-2j}{2}\rfloor}\binom{n-j}{u}\binom{n-j-u}{j+u}(ab)^{j+u}(a+b)^{n-2j-2u}\text{.}
\end{equation}

The Eq.~(\ref{eq:81}) proves the Eq.~(\ref{eq:34}).

\subsection{A Proof of the Eq.~(\ref{eq:16})}

By setting $n:=2n$ in the Eq.~(\ref{eq:81}), it follows that:

\begin{equation}\label{eq:82}
M_B(2n,j,1;a,b)=\binom{2n}{j}\sum_{u=0}^{n-j}\binom{2n-j}{u}\binom{2n-j-u}{j+u}(ab)^{j+u}(a+b)^{2n-2j-2u}\text{.}
\end{equation}

By using the substitution $u=n-j-s$, the Eq.~(\ref{eq:82}) becomes:

\begin{equation}\label{eq:83}
M_B(2n,j,1;a,b)=\sum_{s=0}^{n-j}\binom{2n}{j}\binom{2n-j}{n+s}\binom{n+s}{n-s}(ab)^{n-s}(a+b)^{2s}\text{.}
\end{equation}

It is readily verified that

\begin{equation}\label{eq:84}
\binom{2n}{j}\binom{2n-j}{n+s}=\binom{2n}{n+s}\binom{n-s}{j}\text{.}
\end{equation}

By using the Eq.~(\ref{eq:84}), the Eq.~(\ref{eq:83}) becomes:
\begin{equation}\label{eq:85}
M_B(2n,j,1;a,b)=\sum_{s=0}^{n-j}\binom{2n}{n+s}\binom{n-s}{j}\binom{n+s}{2s}(ab)^{n-s}(a+b)^{2s}\text{.}
\end{equation}

By using the Eq.~(\ref{eq:50}) in the Eq.~(\ref{eq:85}), it follows that 
\begin{equation}\label{eq:86}
M_B(2n,j,1;a,b)=\sum_{s=0}^{n-j}[\frac{\binom{2n}{n}\binom{n}{s}^2}{\binom{2s}{s}}](ab)^{n-s}(a+b)^{2s} \binom{n-s}{j}\text{.}
\end{equation}

Let $p$ be a prime number such that divides $a+b$. Due to the fact that $ab$ and $a+b$ are relatively prime, it follows that $p$ does not divide $ab$. Let $j$ be a fixed non-negative integer such that $j \leq n$.

Obviously,
\begin{equation}\label{eq:87}
v_p([\frac{\binom{2n}{n}\binom{n}{s}^2}{\binom{2s}{s}}](ab)^{n-s}(a+b)^{2s} \binom{n-s}{j} )\geq v_p([\frac{\binom{2n}{n}\binom{n}{s}^2}{\binom{2s}{s}}](ab)^{n-s}(a+b)^{2s})
\end{equation}

By using Eqns.~(\ref{eq:52}) and (\ref{eq:57}), we know that
\begin{equation}\label{eq:88}
v_p([\frac{\binom{2n}{n}\binom{n}{s}^2}{\binom{2s}{s}}](ab)^{n-s}(a+b)^{2s})\geq v_p(\binom{2n}{n})\text{;}
\end{equation}
where $0 \leq s \leq n-j$.

By using Eqns.~(\ref{eq:87}) and (\ref{eq:88}), it follows that:
\begin{equation}\label{eq:89}
v_p([\frac{\binom{2n}{n}\binom{n}{s}^2}{\binom{2s}{s}}](ab)^{n-s}(a+b)^{2s}\binom{n-s}{j})\geq v_p(\binom{2n}{n})\text{;}
\end{equation}
where $0 \leq s \leq n-j$.

By using Eqns.~(\ref{eq:86}) and (\ref{eq:89}), it follows that:

\begin{equation}\label{eq:90}
v_p(M_B(2n,j,1;a,b))\geq v_p(\binom{2n}{n})\text{;}
\end{equation}
where $0\leq j \leq n$.

Now, by using Relation \ref{eq:13} and the Eq.~(\ref{eq:90}), it follows that:

\begin{equation}\notag
\omega_{a+b}(M_B(2n,j,1;a,b))=\min (\lfloor \frac{v_p((M_B(2n,j,1;a,b))}{v_p(a+b)} \rfloor: p \text{ is a prime such that divides } a+b )
\end{equation}
\begin{align}
&\geq\min (\lfloor \frac{v_p(\binom{2n}{n})}{v_p(a+b)} \rfloor: p \text{ is a prime such that divides } a+b )\label{eq:91}
\end{align}

By using Eqns.~(\ref{eq:13}) and (\ref{eq:91}), it follows that:
\begin{equation}
\omega_{a+b}(M_B(2n,j,1;a,b))\geq \omega_{a+b}{(\binom{2n}{n})}\text{;}\label{eq:92}
\end{equation}
where $0\leq j \leq n$.

By the Eq.~(\ref{eq:92}), it follows that:

\begin{equation}\label{eq:93}
(a+b)^{\omega_{a+b}{(\binom{2n}{n})}} \text{ divides } M_B(2n,j,1;a,b)\text{;}
\end{equation}
for any $j$ such that $0\leq j \leq n$.

Note that $q(2n,a,b)=(a+b)^{\omega_{a+b}{(\binom{2n}{n})}}$. By  Corollary \ref{cor:6}, it follows that:
\begin{equation}\label{eq:94}
(a+b)^{\omega_{a+b}{(\binom{2n}{n})}} \text{ divides } M_B(2n,j,t;a,b)\text{;}
\end{equation}
for any $j$ and for any natural number $t$ such that $0\leq j \leq n$.

Furthermore, by Corollary \ref{cor:6}, it follows that:

\begin{equation}\label{eq:95}
(a+b)^{\omega_{a+b}{(\binom{2n}{n})}} \text{ divides } B(2n,m,a,b)\text{;}
\end{equation}
for any natural number $m$ such that $ m\geq 2$.

The Eq.~(\ref{eq:95}) proves the Eq.~(\ref{eq:16}).

\subsection{A Proof of the Eq.~(\ref{eq:17})}

By setting $n:=2n+1$ in the Eq.~(\ref{eq:81}), it follows that:

\begin{equation}\label{eq:96}
M_B(2n+1,j,1;a,b)=\binom{2n+1}{j}\sum_{u=0}^{n-j}\binom{2n+1-j}{u}\binom{2n+1-j-u}{j+u}(ab)^{j+u}(a+b)^{2n+1-2j-2u}\text{;}
\end{equation}
where $0 \leq j \leq n$.

By using the substitution $u=n-j-s$, the Eq.~(\ref{eq:96}) becomes, as follows:

\begin{equation}\label{eq:97}
M_B(2n+1,j,1;a,b)=(a+b)\sum_{s=0}^{n-j}\binom{2n+1}{j}\binom{2n+1-j}{n+s+1}\binom{n+1+s}{2s+1}(ab)^{n-s}(a+b)^{2s}\text{.}
\end{equation}

It is readily verified that:

\begin{equation}\label{eq:98}
\binom{2n+1}{j}\binom{2n+1-j}{n+s+1}=\binom{2n+1}{n+s+1}\binom{n-s}{j}\text{.}
\end{equation}

By the Eq.~(\ref{eq:98}), the Eq.~(\ref{eq:97}) becomes, as follows:

\begin{equation}\label{eq:99}
M_B(2n+1,j,1;a,b)=(a+b)\sum_{s=0}^{n-j}\binom{2n+1}{n+s+1}\binom{n+1+s}{2s+1}\binom{n-s}{j} (ab)^{n-s}(a+b)^{2s}\text{.}
\end{equation}

By using the Eq.~(\ref{eq:65}) in the Eq.~(\ref{eq:99}), it follows that $M_B(2n+1,j,1;a,b)$ is equal to the following sum:

\begin{equation}\label{eq:100}
(a+b)\sum_{s=0}^{n-j}[\frac{(2n+1)\binom{2n}{n}\binom{n}{s}^2}{(2s+1)\binom{2s}{s}}]\binom{n-s}{j} (ab)^{n-s}(a+b)^{2s}\text{.}
\end{equation}

Let $P(2n+1,a,b)$ denote the sum $\sum_{s=0}^{n-j}[\frac{(2n+1)\binom{2n}{n}\binom{n}{s}^2}{(2s+1)\binom{2s}{s}}]\binom{n-s}{j} (ab)^{n-s}(a+b)^{2s}$.

By the Eq.~(\ref{eq:100}) , it follows that:

\begin{equation}\label{eq:101}
\omega_{a+b}(M_B(2n+1,j,1;a,b))=1+\omega_{a+b}(P(2n+1,a,b))\text{.}
\end{equation}

Let $p$ be a prime number such that $p$ divides $a+b$.

Obviously,

\begin{equation}\label{eq:102}
v_p([\frac{(2n+1)\binom{2n}{n}\binom{n}{s}^2}{(2s+1)\binom{2s}{s}}]\binom{n-s}{j} (ab)^{n-s}(a+b)^{2s}) \geq v_p([\frac{(2n+1)\binom{2n}{n}\binom{n}{s}^2}{(2s+1)\binom{2s}{s}}](ab)^{n-s}(a+b)^{2s})\text{.}
\end{equation}

By using Eqns.~(\ref{eq:69}) and (\ref{eq:73}), we know that:

\begin{equation}\label{eq:103}
v_p([ \frac{(2n+1)\binom{2n}{n}\binom{n}{s}^2}{(2s+1)\binom{2s}{s}}] (a+b)^{2s}\cdot (ab)^{n-s} )\geq  v_p((2n+1)\binom{2n}{n})\text{;}
\end{equation}
where $0\leq s \leq n$.

By Eqns.~(\ref{eq:102}) and (\ref{eq:103}), it follows that:

\begin{equation}\label{eq:104}
v_p([\frac{(2n+1)\binom{2n}{n}\binom{n}{s}^2}{(2s+1)\binom{2s}{s}}]\binom{n-s}{j} (ab)^{n-s}(a+b)^{2s}) \geq  v_p((2n+1)\binom{2n}{n})\text{;}
\end{equation}
where $0\leq s \leq n$.

Then, it follows that:

\begin{equation}\label{eq:105}
v_p(P(2n+1,a,b))\geq  v_p((2n+1)\binom{2n}{n})\text{.}
\end{equation}

Now, by using Relation \ref{eq:13} and the Eq.~(\ref{eq:105}), it follows that:

\begin{equation}
\omega_{a+b}(P(2n+1,a,b))=\min (\lfloor \frac{v_p((P(2n+1,a,b))}{v_p(a+b)} \rfloor: p \text{ is a prime such that  divides } a+b )\notag
\end{equation}
\begin{equation}
\geq\min (\lfloor \frac{v_p((2n+1)\binom{2n}{n})}{v_p(a+b)} \rfloor: p \text{ is a prime such that  divides } a+b )\text{.}\label{eq:106}
\end{equation}

By using Relation \ref{eq:13} and the Eq.~(\ref{eq:106}), it follows that:

\begin{equation}\label{eq:107}
\omega_{a+b}(P(2n+1,a,b))\geq\omega_{a+b}((2n+1)\binom{2n}{n})\text{.}
\end{equation}

By Eqns.~(\ref{eq:101}) and (\ref{eq:107}), it follows that:

\begin{equation}\label{eq:108}
\omega_{a+b}(M_B(2n+1,j,1;a,b))\geq 1+\omega_{a+b}((2n+1)\binom{2n}{n})\text{.}
\end{equation}

By the Eq.~(\ref{eq:108}), it follows that
\begin{equation}\label{eq:109}
(a+b)^{1+\omega_{a+b}{((2n+1)\binom{2n}{n})}} \text{ divides } M_B(2n+1,j,1;a,b)\text{;}
\end{equation}
for any $j$ such that $0\leq j \leq n$.

Note that $q(2n+1,a,b)=(a+b)^{1+\omega_{a+b}{((2n+1)\binom{2n}{n})}}$. By  Corollary \ref{cor:6}, it follows that:
\begin{equation}\label{eq:110}
(a+b)^{1+\omega_{a+b}{((2n+1)\binom{2n}{n})}} \text{ divides } M_B(2n+1,j,t;a,b)\text{;}
\end{equation}
for any $j$ and for any natural number $t$ such that $0\leq j \leq n$.

Furthermore, by Corollary \ref{cor:6}, it follows that:

\begin{equation}\label{eq:111}
(a+b)^{1+\omega_{a+b}{((2n+1)\binom{2n}{n})}} \text{ divides } B(2n+1,m,a,b)\text{;}
\end{equation}
for any natural number $m$ such that $ m\geq 2$.

The Eq.~(\ref{eq:111}) proves the Eq.~(\ref{eq:17}). This completes the proof of Theorem \ref{t:2}.

\section{A Proof of Corolary \ref{cor:1}}\label{s:7}

It is well-known that:
\begin{equation}\label{eq:112}
\sum_{k=0}^{n}\binom{n}{k}^2=\binom{2n}{n}\text{.}
\end{equation}

By the Eq.~(\ref{eq:112}), it follows that:

\begin{equation}\label{eq:113}
B(n,2,1,1)=\binom{2n}{n}\text{.}
\end{equation}

By setting $n:=2n$, $a:=1$, and $b:=1$ in the Eq.~(\ref{eq:14}) of Theorem \ref{t:1}, it follows that:
\begin{equation}\label{eq:114}
v_2(B(2n,2,1,1))=v_2(\binom{2n}{n})\text{.}
\end{equation}

Note that we use notation $v_2$ instead of $\omega_2$ because $2$ is a prime number.

By the Eq.~(\ref{eq:113}), it follows that:

\begin{equation}\label{eq:115}
B(2n,2,1,1)=\binom{4n}{2n}\text{.}
\end{equation}

By using the Eq.~(\ref{eq:115}), the Eq.~(\ref{eq:114}) becomes

\begin{equation}\notag
v_2(\binom{4n}{2n})=v_2(\binom{2n}{n})\text{.}
\end{equation}

The last equation above proves the Eq.~(\ref{eq:19}).

Furthermore, by setting $n:=2n+1$, $a:=1$, and $b:=1$  in the Eq.~(\ref{eq:15}) of Theorem \ref{t:1}, it follows that:
\begin{equation}\label{eq:116}
v_2(B(2n+1,2,1,1))=1+v_2((2n+1)\binom{2n}{n})\text{.}
\end{equation}

By the Eq.~(\ref{eq:113}), it follows that:

\begin{equation}\label{eq:117}
B(2n+1,2,1,1)=\binom{4n+2}{2n+1}\text{.}
\end{equation}

By using the Eq.~(\ref{eq:117}), the Eq.~(\ref{eq:116}) becomes, as follows:

\begin{align*}
v_2(\binom{4n+2}{2n+1})&=1+v_2((2n+1)\binom{2n}{n})\text{,}\\
&=1+v_2(2n+1)+v_2(\binom{2n}{n})\text{,}\\
&=1+v_2(\binom{2n}{n})\text{.}
\end{align*}

The last equation above proves the Eq.~(\ref{eq:18}).
This completes the proof of Corollary \ref{cor:1}.

\section{A Proof of Corollary \ref{cor:2}}\label{s:8}

For the proof of Corollary \ref{cor:2}, we use Theorem \ref{t:2}.
It is known that:
\begin{equation}\notag
f_n=B(n,3,1,1)\text{;}
\end{equation}
where $f_n$ is the $n$-th Franel number.

By setting $n:=2n$, $m:=3$, $a=b=1$ in the Eq.~(\ref{eq:16}), 
it follows that

\begin{equation}\notag
v_2(f_{2n})\geq v_2(\binom{2n}{n})\text{.}
\end{equation}

This proves the Eq.~(\ref{eq:21}).

By setting $n:=2n+1$, $m:=3$, $a=b=1$ in the Eq.~(\ref{eq:17}), 
it follows that

\begin{equation}\notag
v_2(f_{2n+1})\geq 1+ v_2((2n+1)\binom{2n}{n})\text{.}
\end{equation}

Obviously,

\begin{equation}\notag
v_2((2n+1)\binom{2n}{n})= v_2(\binom{2n}{n})\text{.}
\end{equation}

By using two last equations above, it follows that:

\begin{equation}\notag
v_2(f_{2n+1})\geq 1+ v_2(\binom{2n}{n})\text{.}
\end{equation}

This proves the Eq.~(\ref{eq:20}).

\section{ A Proof of Theorem \ref{t:3}}\label{s:9}
It is known that:
\begin{equation}\notag
D_n=B(n,2,1,2)\text{;}
\end{equation}
where $D_n$ is  $n$-th central Dellanoy number.

By setting $n:=2n+1$, $a:=1$, and $b:=2$ in the Eq.~(\ref{eq:15}), it follows that:
\begin{align*}
v_3(D_{2n+1})&=1+v_3((2n+1)\binom{2n}{n})\text{,}\\
&=1+v_3((2n+1))+v_3(\binom{2n}{n})\text{.}
\end{align*}

The last equation above proves the Eq.~(\ref{eq:22}).

Note that we use notation $v_3$ instead of $\omega_3$ because $3$ is a prime number.

By setting $n:=2n$, $a:=1$, and $b:=2$ in the Eq.~(\ref{eq:14}), it follows that:
\begin{equation}\notag
v_3(D_{2n})=v_3(\binom{2n}{n})\text{.}
\end{equation}

The last equation above proves the Eq.~(\ref{eq:23}).

This completes the proof of Theorem \ref{t:3}.

\section{ Proofs of Theorem \ref{t:5} and Corollary \ref{cor:3}}\label{s:10}

Proof of Theorem \ref{t:5} directly follows from the our proof of Theorem \ref{t:1}.
Just replace $ab$ with $a$, and $a+b$ with $b$ in the proof of Theorem \ref{t:1}.

Also Corollary \ref{cor:3} follows from Theorem \ref{t:1} and the Eq.~(\ref{eq:3}).

\section{ A Proof of Theorem \ref{t:6}}\label{s:11}

Proof of Theorem \ref{t:6} is similar to the proof of Theorem \ref{t:1}.

Note that:

\begin{equation}\label{eq:118}
M_n(a,b)=\sum_{k=0}^{\lfloor \frac{n}{2}\rfloor}\binom{n}{k}\binom{n-k}{k}\cdot \frac{1}{k+1} \cdot a^k\cdot b^{n-2k}\text.{}
\end{equation}

\subsection{A Proof of the Eq.~(\ref{eq:28.3})}\label{s:11.1}

By setting $n:=2n$ in the Eq.~(\ref{eq:118}), it follows that:

\begin{equation}\label{eq:119}
M_{2n}(a,b)=\sum_{k=0}^{n}\binom{2n}{k}\binom{2n-k}{k}\cdot \frac{1}{k+1} \cdot a^k\cdot b^{2n-2k}\text{.}
\end{equation}

By using the substitution $k=n-s$, the Eq.~(\ref{eq:119}) becomes gradually:

\begin{align}
M_{2n}(a,b)&=\sum_{s=0}^{n}\binom{2n}{n-s}\binom{n+s}{n-s}\cdot \frac{1}{n+1-s}\cdot b^{2s}\cdot  a^{n-s}\text{,}\notag\\
&=\sum_{s=0}^{n}\binom{2n}{n+s}\binom{n+s}{2s}\cdot \frac{1}{n+1-s}\cdot  b^{2s}\cdot a^{n-s}\text{.}\label{eq:120}
\end{align}

By using the Eq.~(\ref{eq:50}), the Eq.~(\ref{eq:120}) becomes, as follows:

\begin{equation}\label{eq:121}
M_{2n}(a,b)=\sum_{s=0}^{n}\bigl{(} \frac{\binom{2n}{n}\binom{n}{s}^2}{\binom{2s}{s}} \bigr{)}\cdot \frac{1}{n+1-s}\cdot b^{2s}\cdot a^{n-s}\text{.}
\end{equation}

It is readily verified that:

\begin{equation}\label{eq:122}
\binom{n}{s}\cdot \frac{1}{n+1-s}=\frac{1}{n+1}\cdot \binom{n+1}{s}\text{.}
\end{equation}

By using the  Eq.~(\ref{eq:122}), the Eq.~(\ref{eq:121}) becomes, as follows:

\begin{equation}\label{eq:123}
M_{2n}(a,b)=\sum_{s=0}^{n}\bigl{(}\frac{ C_n\binom{n}{s}\binom{n+1}{s}}{\binom{2s}{s}} \bigr{)}\cdot b^{2s}\cdot a^{n-s}\text{.}
\end{equation}

Let $p$ be a prime number such that divides $b$.

For $s=0$, the first summand in the right-side of the Eq.~(\ref{eq:123}) is $C_n \cdot a^n$.  Since $a$ and $b$ are relatively prime, it follows that:

\begin{equation}\label{eq:124}
v_p(C_n \cdot a^n)=v_p(C_n)\text{.}
\end{equation}

Let $ 1\leq s \leq n$. Then, for $s+1$-th summand  in the right-side of the Eq.~(\ref{eq:123}), we have:

\begin{equation}\label{eq:125}
v_p(\bigl{(}\frac{ C_n\binom{n}{s}\binom{n+1}{s}}{\binom{2s}{s}} \bigr{)}\cdot b^{2s}\cdot a^{n-s})\geq v_p(C_n)-v_p(\binom{2s}{s})+2s\cdot v_p(b)\text{.}
\end{equation}

By using the Ineq.~(\ref{eq:53}), from the Eq.~(\ref{eq:125}), we obtain that:

\begin{equation}\label{eq:126}
v_p(\bigl{(}\frac{ C_n\binom{n}{s}\binom{n+1}{s}}{\binom{2s}{s}} \bigr{)}\cdot b^{2s}\cdot a^{n-s}) \geq v_p(C_n)+s\cdot(2v_p(b)-1)\text{.}
\end{equation}

Obviously, for $1\leq s \leq n$, 

\begin{equation}\label{eq:127}
v_p(C_n)+s\cdot(2v_p(b)-1)\geq v_p(C_n)+s \geq  v_p(C_n)+1 \text{.}
\end{equation}

By Ineq.~(\ref{eq:126}) and (\ref{eq:127}), for $1\leq s \leq n$, it follows that:
\begin{equation}\label{eq:128}
v_p(\bigl{(}\frac{ C_n\binom{n}{s}\binom{n+1}{s}}{\binom{2s}{s}} \bigr{)}\cdot b^{2s}\cdot a^{n-s})\geq  v_p(C_n)+1\text{.}
\end{equation}

By using Eqns. (\ref{eq:124}) and (\ref{eq:128}), we obtain that:

\begin{equation}\label{eq:129}
v_p(M_{2n}(a,b))=v_p(C_n)\text{,}
\end{equation}
where $p$ is a prime number that divides $b$.

Furthermore, by using Relation \ref{eq:13}, it follows that:

\begin{align}
\omega_{b}(M_{2n}(a,b))&=\min (\lfloor \frac{v_p(M_{2n}(a,b))}{v_p(b)} \rfloor: p \text{ is a prime such that divides } b )\label{eq:130}\\
&=\min (\lfloor \frac{v_p(C_n)}{v_p(b)} \rfloor: p \text{ is a prime such that divides } b )\label{eq:131}\\
&=\omega_{b}{(C_n)}\text{.}\label{eq:132}
\end{align}

The Eq.~(\ref{eq:132}) completes the proof of the Eq.~(\ref{eq:28.3}).

\subsection{A Proof of the Eq.~(\ref{eq:28.4})}\label{s:11.2}

By setting $n:=2n+1$ in the Eq.~(\ref{eq:118}), we obtain that:

\begin{align}\
M_{2n+1}(a,b)&=\sum_{k=0}^{n}\binom{2n+1}{k}\binom{2n+1-k}{k}\cdot \frac{1}{k+1}\cdot a^k\cdot b^{2n+1-2k}\text{,}\notag\\
&=b\sum_{k=0}^{n}\binom{2n+1}{k}\binom{2n+1-k}{k}\cdot \frac{1}{k+1}\cdot a^k\cdot b^{2n-2k}\text{.}\label{eq:133}
\end{align}

By using the substitution $k=n-s$, the Eq.~(\ref{eq:133}) becomes gradually:

\begin{align}
M_{2n+1}(a,b)&=b\sum_{s=0}^{n}\binom{2n+1}{n-s}\binom{n+s+1}{n-s}\cdot \frac{1}{n-s+1}\cdot b^{2s}\cdot  a^{n-s}\text{,}\notag\\
&=b\sum_{s=0}^{n}\binom{2n+1}{n+s+1}\binom{n+s+1}{2s+1} \frac{1}{n+1-s}\cdot b^{2s}\cdot a^{n-s}\text{.}\label{eq:134}
\end{align}

By using the Eq.~(\ref{eq:65}), the Eq.~(\ref{eq:134}) becomes, as follows:

\begin{equation}\label{eq:135}
M_{2n+1}(a,b)=b\sum_{s=0}^{n}\bigl{[} \frac{(2n+1)\binom{2n}{n}\binom{n}{s}^2}{(2s+1)\binom{2s}{s}} \bigr{]}\cdot \frac{1}{n+1-s}\cdot b^{2s}\cdot a^{n-s}\text{.}
\end{equation}

By using the Eq.~(\ref{eq:122}), it follows that

\begin{equation}\label{eq:136}
M_{2n+1}(a,b)=b\sum_{s=0}^{n}\bigl{[} \frac{(2n+1)\cdot C_n\binom{n}{s}\binom{n+1}{s}}{(2s+1)\binom{2s}{s}} \bigr{]}\cdot  b^{2s}\cdot a^{n-s}\text{.}
\end{equation}

The rest of the proof is similar to the proof of  the Eq.~(\ref{eq:15}) of Theorem \ref{t:1}.

Let $p$ be a prime that divides $b$. Then, by using Lemma \ref{l:1}, it can be shown that:

\begin{equation}\label{eq:137}
v_p(\sum_{s=0}^{n}\bigl{[} \frac{(2n+1)\cdot C_n\binom{n}{s}\binom{n+1}{s}}{(2s+1)\binom{2s}{s}} \bigr{]}\cdot  b^{2s}\cdot a^{n-s})=v_p((2n+1)\cdot C_n)\text{.}
\end{equation}

This completes the proof of Theorem \ref{t:6}.

\section{ A Proof of Theorem \ref{t:4}}\label{s:12}

We use Theorem \ref{t:6} and the Eq.~(\ref{eq:8.4}).

By setting $a:=2$ and $b:=3$ in the Eq.~(\ref{eq:28.3}), it follows that:

\begin{equation}\label{eq:138}
v_3(M_{2n}(2,3))=v_3(C_n)\text{.}
\end{equation}

By the Eq.~(\ref{eq:8.4}), it follows that

\begin{equation}\label{eq:139}
s_{2n+1}=M_{2n}(2,3)\text{.}
\end{equation}

By using Eqns.~(\ref{eq:138}) and (\ref{eq:139}), it follows that:

\begin{equation}\label{eq:140}
v_3(s_{2n+1})=v_3(C_n)\text{,}
\end{equation}
where $n$ is a non-negative integer.

The Eq.~(\ref{eq:24}) follows from the Eq.~(\ref{eq:140}) and the fact that $v_3(S_n)=v_3(s_n)$.

Similarly, by setting $a:=2$ and $b:=3$ in the Eq.~(\ref{eq:28.4}), it follows that:

\begin{equation}\label{eq:141}
v_3(M_{2n+1}(2,3))=1+v_3((2n+1)\cdot C_n)\text{.}
\end{equation}

By the Eq.~(\ref{eq:8.4}), it follows that

\begin{equation}\label{eq:142}
s_{2n+2}=M_{2n+1}(2,3)\text{.}
\end{equation}

By using Eqns.~(\ref{eq:141}) and (\ref{eq:142}), it follows that:

\begin{equation}\label{eq:143}
v_3(s_{2n+2})=1+v_3(2n+1)+v_3(C_n)\text{,}
\end{equation}
where $n$ is a non-negative integer.

Now, the Eq.~(\ref{eq:25}) follows from the Eq.~(\ref{eq:143}) and the fact that $v_3(S_n)=v_3(s_n)$.

This completes the proof of Theorem \ref{t:4}.

\begin{remark}\label{r:5}
Theorem \ref{t:4} can be proved by using Theorem \ref{eq:3}, Eqns.~(\ref{eq:5}), (\ref{eq:6}), and Kummer's theorem.
\end{remark}

\section{Concluding Remarks}\label{s:13}

\begin{remark}\label{r:3}
It is known \cite[page 2]{Wang} and \cite[Eq.~(1.5)]{Zhao} that

\begin{equation}\label{eq:144}
H_{n}=M_n(1,3)\text{,}
\end{equation}
where $H_n$  are restricted hexagonals  numbers \cite{Read}.

By using Theorem \ref{t:6}, it can be shown that:

\begin{align}
v_3(H_{2n})&=v_3(C_n)\text{,}\label{eq:145}\\
v_3(H_{2n+1})&=1+v_3(2n+1)+v_3(C_n)\text{.}\label{eq:146}
\end{align}
\end{remark}

\begin{remark}\label{r:4}
It is known \cite[page 2]{Wang} that

\begin{equation}\label{eq:147}
C_{n+1}=M_n(1,2)\text{,}
\end{equation}
where $C_n$  denotes the $n$th Catalan number.

By using Theorem \ref{t:6}, it can be shown that:

\begin{align}
v_2(C_{2n+1})&=v_2(C_n)\text{,}\label{eq:148}\\
v_2(C_{2n+2})&=1+v_2(C_n)\text{.}\label{eq:149}
\end{align}
\end{remark}

\begin{remark}\label{r:6}

By using Theorem \ref{t:1} and Kummer's theorem \cite[Theorem 1.~2, p.\ 2]{Mihet}, we can efficiently calculate $\omega_{a+b}(B(n,2,a,b))$.
Let us consider the following sum: \[B(4046,2,37,62)=\sum_{k=0}^{4046}\binom{4046}{k}^2\cdot 37^{4046-k}\cdot 62^k\text{.}\]
Here $a:=37$, $b:=62$, $2n:=4046$, and $a+b:=99$. Obviously, $\gcd(37, 62)=1$.

By the Eq.~(\ref{eq:14}) of Theorem \ref{t:1}, it follows that:
\begin{equation*}
\omega_{99}(B(4046,2,37,62))=\omega_{99}(\binom{4046}{2023})\text{.}
\end{equation*}
Since $99=3^2\cdot 11$, it follows $v_3(99)=2$ and $v_{11}(99)=1$. We need to find $v_3(\binom{4046}{2023})$ and $v_{11}(\binom{4046}{2023})$.

The expansion of $2023$ in the base $3$ is $(2202221)_3$. If we add the number $(2202221)_3$ to itself, there are five   `carry-overs' in total.  By Kummer's theorem,  $v_3(\binom{4046}{2023})=5$.

The expansion of $2023$ in the base $11$ is $(157A)_{11}$. If we add the number  $(157A)_{11}$. to itself, there are three `carry-overs' in total.  By Kummer's theorem,  $v_{11}(\binom{4046}{2023})=3$.

By the Eq.~(\ref{eq:13}), it follows that:

\begin{align*}
\omega_{99}(\binom{4046}{2023})&=\min (\lfloor \frac{v_{11}(\binom{4046}{2023})}{v_{11}(99)} \rfloor,  \lfloor \frac{v_3(\binom{4046}{2023})}{v_3(99)} \rfloor)\\
&=\min (\lfloor \frac{3}{1} \rfloor,  \lfloor \frac{5}{2} \rfloor)\\
&=\min (3, 2)\\
&=2\text{.}
\end{align*}

By Theorem \ref{t:1}, it follows that 
\begin{equation*}
\omega_{99}(B(4046,2,37,62))=2\text{.}
\end{equation*}

Moreover, by Theorem \ref{t:2}, it follows that the sum $B(4046,m,37,62)$ is always divisible by $99^2$ for all natural numbers $m$ such that $m \geq 2$.

Also, let us consider the following sum: \[B(4047,2,37,62)=\sum_{k=0}^{4047}\binom{4047}{k}^2\cdot 37^{4047-k}\cdot 62^k\text{.}\]
By the Eq.~(\ref{eq:15}) of Theorem \ref{t:1}, it follows that:
\begin{equation*}
\omega_{99}(B(4046,2,37,62))=1+\omega_{99}(4047\cdot \binom{4046}{2023})\text{.}
\end{equation*}

We have:

\begin{align*}
v_3(4047\cdot \binom{4046}{2023})&=v_3(4047)+v_3( \binom{4046}{2023})\\
&=1+5=6\text{.}
\end{align*}

Furthermore, 
\begin{align*}
v_{11}(4047\cdot \binom{4046}{2023})&=v_{11}(4047)+v_{11}(\binom{4046}{2023})\\
&=0+3=3\text{.}
\end{align*}

By the Eq.~(\ref{eq:13}), it follows that:

\begin{align*}
\omega_{99}(4047\binom{4046}{2023})&=\min (\lfloor \frac{v_{11}(4047\binom{4046}{2023})}{v_{11}(99)} \rfloor,  \lfloor \frac{v_3(4047\binom{4046}{2023})}{v_3(99)} \rfloor)\\
&=\min (\lfloor \frac{3}{1} \rfloor,  \lfloor \frac{6}{2} \rfloor)\\
&=\min (3, 3)\\
&=3\text{.}
\end{align*}

By the Eq.~(\ref{eq:15}) of Theorem \ref{t:1}, it follows that:
\begin{align*}
\omega_{99}(B(4046,2,37,62))&=1+\omega_{99}(4047\cdot \binom{4046}{2023})\\
&=1+3=4\text{.}
\end{align*}

It follows that: 
\[ \omega_{99}(B(4046,2,37,62))=4\text{.}\]

Moreover, by Theorem \ref{t:2}, it follows that the sum $B(4047,m,37,62)$ is always divisible by $99^4$ for all natural numbers $m$ such that $m \geq 2$.
\end{remark}

\begin{remark}\label{r:7}
By using Eqns.~(\ref{eq:32}) and (\ref{eq:34}), we can derive a slight generalization of MacMahon's formula (\ref{eq:8}).
It can be shown that:

\begin{equation}\label{eq:150}
M_B(n,j,2;a,b)=\binom{n}{j}\sum_{l=0}^{\lfloor \frac{n}{2}\rfloor}\binom{n}{2l}\binom{2l}{l}\binom{n+l-j}{n}(ab)^{l}(a+b)^{n-2l}\text{.}
\end{equation}
By setting $j:=0$ in the Eq.~(\ref{eq:150}), the Eq.~(\ref{eq:8}) follows.
\end{remark}

\section*{Acknowledgments}

I want to thank to my former teacher Aleksandar Stankov Leko.

\end{document}